\input amstex
\documentstyle{amsppt}
\magnification=1200
\vcorrection{-0.2in}
\NoRunningHeads
\NoBlackBoxes

\define\a{\alpha}
\define\al{\aligned}
\define\eda{\endaligned}

\define\bl{\bigl}
\define\br{\bigr}


\define\de{\delta}

\define\dem{\demo{Proof}}

\define\edo{\enddemo}

\define\fr{\frac}
\define\g{\gamma}

\define\la{\lambda}
\define\lf{\left}

\define\p{\partial}
\define\rg{\right}

\define\={\equiv}
\topmatter
\title Positive quaternionic K\"ahler manifolds and symmetry rank
\endtitle
\author Fuquan Fang \footnote{Supported by NSF Grant 19925104 of
China, 973 project of Foundation Science of China, RFDP}
\endauthor
\address
Nankai Institute of Mathematics, Nankai University,
Tianjin
300071, P.R.C.
\newline
Instituto de Matematica, Universidade Federal Fluminense, Niteroi, R.J.,
Brasil
\endaddress
\email ffang\@nankai.edu.cn ; fuquanfang\@eyou.com
\endemail
\abstract A quaternionic K\"ahler manifold $M$ is called {\it positive} if
it has positive scalar curvature. The main purpose of this paper is to prove
several connectedness theorems for quaternionic immersions in a quaternionic
K\"ahler manifold, e.g. the Barth-Lefschetz type connectedness theorem
for quaternionic submanifolds in a positive quaternionic K\"ahler manifold.
As applications we prove that, among others, a $4m$-dimensional positive
quaternionic K\"ahler manifold with symmetry rank at least $(m-2)$ must be
either isometric to $\Bbb HP^m$ or $Gr_2(\Bbb C^{m+2})$, if $m\ge 10$.
\endabstract
\endtopmatter
\document
\baselineskip=.25in

\vskip2mm

\head 0. Introduction
\endhead

\vskip4mm

A quaternionic K\"ahler manifold $M$ is an oriented Riemannian $4n$-manifold,
$n\ge 2$, whose holonomy group is contained in $Sp(n)Sp(1)\subset SO(4n)$.
If $n=1$ we add the condition that $M$ is Einstein and self dual.
Equivalently, there exists a $3$-dimensional subbundle $S$, of the
endmorphism bundle $\text{End}(TM, TM)$ locally
generated by three anti-commuting almost complex structures $I, J, K=IJ$
so that the Levi-Civita connection preserves $S$. It is well-known [Ber] that a
quaternionic K\"ahler manifold $M$ is always Einstein, and is necessarily
locally hyperK\"ahler if its Ricci tensor vanishes. Motivated by the Penrose
twistor construction (cf. [AHS] [Hi]),  Salamon [Sa] developped the important
twistor space theory for quaternionic K\"ahler manifolds, showing that the
unit sphere bundle of $S$, called the twistor space $Z$, admits a complex
structure so that the fiber of the $\Bbb P^1$-fibration $p: Z\to M$ is a
rational curve. A quaternionic K\"ahler manifold $M$ is called {\it positive }
if it has positive scalar curvature. By [Hi] (for $n=1$) and [Sa]
(for $n\ge 2$, compare [Le] [Le-Sa]) a positive quaternionic K\"ahler
manifold $M$ has twistor space $Z$ a complex Fano manifold.  Hitchin [Hi]
proved a positive quaternionic K\"ahler $4$-manifold $M$ must be isometric to
$\Bbb CP^2$ or $S^4$. Hitchin's work was extended by Poon-Salamon [PS] to
dimension $8$, showing that a positive quaternionic K\"ahler $8$-manifold
$M$ must be isometric to $\Bbb HP^2$, $Gr_2(\Bbb C^{4})$ or $G_2/SO(4)$.
This leads to the Salamon-Lebrun conjecture:

{\it Every positive quaternionic
K\"ahler manifold is a quaternionic symmetric space.}

Very recently, the conjecture was further verified for $n=3$ in [HH], using
approach initiated in [Sa] [PS] (compare [LeSa]).
For a positive quaternionic K\"ahler manifold $M$, Salamon [Sa] proved that
the dimension of its isometry group is equal to the index of certain twisted
Dirac operator, which by the Atiyah-Singer index theorem, is a characteristic
number of $M$ coupled with the Kraines $4$-form $\Omega $ (in analog with the
K\"ahler form), and it was applied to prove the isometry group of $M$ is
large in lower dimensions (up to dimension $16$).

By [LeSa] any positive quaternionic K\"ahler $4n$-manifold $M$ is simply
connected and the second homotopy group $\pi _2(M)$ is a finite group or
$\Bbb Z$, and $M$ is isometric to $\Bbb HP^n$ or $Gr_2(\Bbb C^{n+2})$
according to $\pi _2(M)=0$ or $\Bbb Z$.  

The main purpose of this paper is to prove several connectedness theorems
for positive quaternionic submanifolds, by using Morse theory on
path spaces, e.g. the analog of the well-known Lefschetz hyperplane section
theorem [Le] (cf. [FL][Fu]), Barth-Lefschetz connectedness theorem,
among others. In particular, our results show that a quaternionic K\"ahler
submanifold of small codimension in a positive quaternionic K\"ahler manifold
shares the homotopy groups up to a certain range.

By Gray [Gr] any quaternionic K\"ahler submanifold must be totally geodesic,
which often shows up as fixed point components of isometric actions. Our
connectedness theorems may be applied to study positive quaternionic K\"ahler
manifold in terms of informations on its isometry group. This approach
dates back to the work [PS] for $n=2$ [HH] for $n=3$ to proving the action is
transitive, and [DS] [PV] for cohomogeneity one actions (and hence the
isometry group must be very large). [Bi] classified positive quaternionic
K\"ahler $4n$-manifolds with isometry rank $n+1$, using an approach on
hyper-K\"ahler quantizations. Based on our connectedness theorems
we will prove that, using a much more direct and transparent approach,
a positive quaternionic K\"ahler $4n$-manifolds with isometry rank
$\ge n-2$ must be either isometric to $\Bbb HP^n$ or $Gr_2(\Bbb C^{n+2})$,
if $n\ge 10$.
                
For the sake of simplicity manifolds, $M, N_1, N_2$ in the paper,
are all closed and connected manifolds. Let $M, N$ be quaternionic
manifolds. An immersion $f: N\to M$ is called {\it quaternionic} if $f$
preserves the quaternionic structures.

\vskip2mm

\proclaim{Theorem A}

Let $M$ be a positive quaternionic K\"ahler manifold of dimension $4m$.
Assume $f=(f_1,f_2): N\to M\times M$, where $N=N_1\times N_2$ and
$f_i: N_i\to M$ are quaternionic immersions of compact quaternionic K\"ahler
manifolds of dimensions $4n_i$, $i=1,2$. Let $\Delta $ be the diagonal of
$M\times M$. Set $n=n_1+n_2$. Then:

\noindent (A1) If $n\ge m$, then $f^{-1}(\Delta)$ is nonempty.

\noindent (A2) If $n\ge m+1$, then $f^{-1}(\Delta)$ is connected.

\noindent (A3) If $f$ is an embedding, then for $i\le n-m$
there is a natural isomorphism, $\pi_i(N_1,N_1\cap N_2)\to
\pi_i(M,N_2)$  and a surjection for $i=n-m+1$.
\endproclaim

\vskip2mm

Some comments on Theorem A are in order.

First, (A1) implies immediately the following intersection theorem, which is
due to Marchiafava [Ma] for embedded quaternionic K\"ahler submanifolds.

\vskip2mm

\proclaim{Theorem 0.1}

Let $M$ be a positive quaternionic K\"ahler manifold of dimension $4m$.
Let $f_i:N_i\to M,\ i=1,2,$ be two immersions of quaternionic K\"ahler
manifolds. If $\dim(N_1)+\dim(N_2)\ge 4m$, then
$f_1(N_1)\cap f_2(N_2)\not=\emptyset$.
\endproclaim

\vskip2mm

Secondly, (A2) implies the following regularity result.

\vskip2mm

\proclaim{Theorem 0.2 (Regularity)}

Let $M$ be a positive quaternionic K\"ahler manifold of dimension $4m$.
If $f: N\to M$ be a quaternionic immersion of a quaternionic K\"ahler
manifold of dimension $4n$. Then $f$ is an embedding, provided $2n\ge m+1$.
\endproclaim

\vskip2mm

Theorem 0.2 may be considered as a quaternionic analog of the well-known
Fulton-Hansen immersion theorem ([FH]) for local complete intersections
in $\Bbb CP^m$, which asserts that a finite unramified morphism from an
algebraic variety of dimension $n$ to $\Bbb CP^m$ must be an embedding,
if the dimension satisfies $2n>m$.

Thirdly, (A3) is exactly the analog of the classical Barth-Lefschetz
hyperplane section Theorems. A particular case implies that the analog
of the Lefschetz theorem, namely, the inclusion $i: N\to M$ is
$(2n-m+1)$-connected where $N$ is a quaternionic K\"ahler submanifold of
dimension $4n$.

\vskip2mm

(A3) is sharp, since $\text{Gr}_2(\Bbb C^4)=\widetilde{\text{Gr}}_4
(\Bbb R^6)\subset \widetilde{\text{Gr}}_4(\Bbb R^7)$ is a quaternionic
K\"ahler submanifold, for which the inclusion is $2$-connected but
not $3$-connected, because the Betti numbers $b_2(\text{Gr}_2(\Bbb C^4))=1$,
and $b_2(\widetilde{\text{Gr}}_4(\Bbb R^7))=0$. In general, the natural
inclusion $\widetilde {\text{Gr}}_4(\Bbb R^m)\subset
\widetilde {\text{Gr}}_4 (\Bbb R^{m+1})$ is a quaternionic K\"ahler
submanifold, for which (A3) can not be improved.

\vskip 2mm

For a Riemannian manifold, by the {\it symmetry rank} we mean the rank of its
isometry group. Bielawski [Bi] proved that a positive quaternionic K\"ahler
$4m$-manifold of symmetry rank at least $m+1$ is isometric to
$\Bbb HP^m$ or $Gr_2(\Bbb C^{m+2})$. Now we state a rigidity theorem of
positive quaternionic K\"ahler manifolds in terms its symmetry rank using
(A3), which much enhances Bielawski's theorem. 

\vskip 2mm

\proclaim{Theorem B}

Let $M$ be a positive quaternionic K\"ahler manifold of dimension $4m$.
Then the isometry group $\text{Isom}(M)$ has rank (denoted by
$\text{rank}(M)$) at most $(m+1)$, and $M$ is isometric to $\Bbb HP^m$ or
$Gr_2(\Bbb C^{m+2})$ if $\text{rank}(M) \ge m-2$ and $m\ge 10$.
\endproclaim

\vskip2mm

The idea to of proving Theorem B is roughly as follows:  First note that
the $T^r$-action on $M$ must have non-empty fixed point set since the Euler
characteristic $\chi (M)>0$ by [Sa]. Consider the isotropy
representation of $T^r$ at a fixed point $x\in M$, which must be a
representation through the local linear holonomy $Sp(n)Sp(1)$ at $T_xM\cong
\Bbb H^n$. Under the condition on $r$ we will prove that there are a
chain of quaternionic K\"ahler submanifolds in $M$, 
$$M^4\subset M^8\subset M^{12}\cdots \subset M$$
such that every nearby two manifolds have relatively small dimension
difference, and the $12$-dimensional quaternionic K\"ahler manifold $M^{12}$
has either symmetry rank $4$ or an effective $T^3$-action of
{\it quaternionic type} (see Section 2 for the definition). Combining [HH] we
see that $M^{12}=\text{Gr}_2(\Bbb C^5)$ or $\Bbb HP^3$. By Theorem (A3) we
know that $\pi _2(M)$ must be isomorphic to $\pi _2(M^{12})$ through the
chain of quaternionic K\"ahler submanifolds. By now we apply [LeSa] to
conclude the desired result.

\vskip2mm

This paper was partially inspired by the work [FMR], where a connectedness
principle was developped for minimal submanifolds in a positively curved
manifold.

\vskip2mm

The rest of the paper is organized as follows:

\noindent In Section 1, we recall some preliminary about
quaternionic K\"ahler manifolds.

\noindent In Section 2, we prove Theorem B assuming Theorem A.

\noindent In Section 3, we apply Morse theory to suitable path spaces.

\noindent In Section 4, we prove some indices Theorems
for energy functions.

\noindent In Section 5, we prove Theorem A and its corollaries.

\vskip10mm

\noindent {\bf Acknowledgment:}
The author would like to thank S. Salamon for some comment on quaternionic
K\"ahler manifolds, and to thank Universidade Federal Fluminense (UFF)
for the hospitality during the preparation of this work.

\vskip10mm

\vskip10mm

\head 1. Quaternionic K\"ahler manifolds
\endhead

\vskip4mm

In this section we recall some basic results on quaternionic K\"ahler
manifolds needed in later sections.

\vskip 2mm

Let $(M,g)$ be a connected Riemannian $4n$-manifold, $n\ge 2$.
We say that $(M, g)$ is a quaternionic K\"ahler manifold if the holonomy
group is conjugate to a subgroup in $Sp(n)Sp(1)=Sp(n)\times _{\Bbb Z_2}
Sp(1)\subset SO(4n)$.  In dimension $4$ we say $(M, g)$ is
quaternionic K\"ahler manifold if it is Einstein, with non-zero scalar
curvature, and self dual. A quaternionic K\"ahler manifold is {\it positive}
if it has positive scalar curvature.

The quaternionic projective spaces $\Bbb HP^n$, the complex Grassmannian
$Gr_2(\Bbb C^{n+2})$, and the oriented real Grassmannian
$\widetilde {\text{Gr}}_4(\Bbb R^{n+4})$ are quaternionic K\"ahler manifolds
with positive scalar curvature. By Wolf [Wo]
there is exactly one quaternionic symmetric space for each compact simple
Lie algebra, and by Alekseevskii [Al] every compact quaternion homogeneous
space is a symmetric space.

As we mentioned in the introduction, by far quaternionic symmetric spaces are
the only known examples of positive quaternionic K\"ahler manifold.

\vskip 2mm

\proclaim{Theorem 1.1 ([Le-Sa])}

(i) (Fininteness) For any $n$, there are, modulo isometries and rescalings,
only finitely many positive quaternionic K\"ahler $4n$-manifolds.

(ii) (Strong rigidity) Let $(M, g)$ be a positive quaternionic K\"ahler
$4n$-manifold. Then $M$ is simply connected and
$$\pi _2(M)=\cases 0, \text{   } (M,g)=\Bbb HP^n \\
\Bbb Z, \text{  } (M,g)=\text{Gr}_2(\Bbb C^{n+2})\\
\text{finite with $2$-torsion}, \text{  } \text{otherwise}
\endcases$$

\endproclaim

\vskip2mm

The following curvature relations due to Berger [Be2] will be very useful.

\vskip 2mm

\proclaim{Proposition 1.2 ([Ber])}
Let $(M, g)$ be a quaternionic K\"ahler $4n$-manifold with Einstein constant
$\lambda$. Let $I, J, K$ be the locally defined three anticommuting almost
complex structures. Then the curvature operator satisfies the following
identities:
$$\langle  R(X,IX)X, IX\rangle+ \langle  R(X,JX)X,JX\rangle+
\langle  R(X,KX)X,KX\rangle=\frac {3\lambda}{n+2}\| X\|^4$$
$$\langle R(X,Y)X, Y\rangle+\langle R(X,IY)X,IY\rangle+\langle R(X,JY)X,JY
\rangle+ \langle R(X,KY)X,KY\rangle= $$
$=\frac {\lambda}{n+2} \| X\|^2\|Y\|^2$
\endproclaim

\vskip2mm

Proposition 1.2 was used to prove a quaternionic K\"ahler manifold is
Einstein [Ber].

\vskip 2mm

\proclaim{Theorem 1.3 ([Fuj] [NT])}
Let $M$ be a positive quaternionic K\"ahler manifold. Then
the intersection form of $M$ is positive definite.
\endproclaim

\vskip 2mm

A submanifold $N$ in a quaternionic K\"ahler manifold is called a {\it
quaternionic submanifold} if the quaternionic structure (i.e. locally
defined $I, J, K$) preserves the tangent bundle of $N$. 

\vskip 2mm

\proclaim{Proposition 1.4 ([Gr])}
Any quaternionic submanifold in a quaternionic K\"ahler manifold is totally
geodesic and quaternionic K\"ahlerian.
\endproclaim

\vskip2mm

\vskip8mm

\head 2. Symmetry rank of Positive Quaternionic K\"ahler manifold
\endhead

\vskip4mm

Let $M$ be a positive quaternionic K\"ahler manifold of dimension $4n$. We
call the rank of the isometry group $\text{Isom}(M)$ the symmetry rank of $M$,
denoted by $\text{rank}(M)$.  Bielawski [Bi]
obtained that a positive quaternionic K\"ahler manifold of symmetry rank
$n+1$ is isometric to  $\Bbb HP^n$ or $Gr_2(\Bbb C^{n+2})$, via studying
hyperk\"ahler $4n+4$-manifolds with tri-Hamiltonian $\Bbb R^{n+1}$-action.
In this section we prove Theorems B and C assuming Theorem A.
Our proof of Theorem B also give a more transparent and direct proof of
his result.

To start let us note that
$
\text{rank}(M) =n+1$ if $M=\Bbb HP^n$ or $\text{Gr}_2(\Bbb C^{n+2})$, and
$\text{rank}(M) =[\frac n2]+2$ if $M=\widetilde{\text{Gr}}_4(\Bbb R^{n+4})$.

\vskip 4mm

\subhead a). A upper bound for symmetry rank
\endsubhead

\vskip 2mm

Let $G\subset \text{Isom}(M)$ be a connected Lie group. For any
$x\in M$, the isotropy subgroup $G_x$ is a subgroup of the holonomy
group $(Sp(n)Sp(1))_x$ at $x$, by the Kostant's Theorem (compare [PS]
page 366). Therefore the isotropy representation is constituted from the two
homomorphisms $$\rho_x: G_x\to Sp(n)_x, \hskip4mm \bar \rho_x: G_x\to
Sp(1)_x$$
Identifying $T_xM\cong \Bbb H^n$ with a quaternionic structure from the
right, observe that $Sp(n)_x$ (resp. $Sp(1)_x$) acts on $T_xM$ from
the left quaternionically linearly (resp. the right by quaternion
multiplication.)

\vskip 2mm

\proclaim{Lemma 2.1}
Let $M$ be a positive quaternionic K\"ahler manifold of dimension $4n$. Then
$\text{rank}(M)\le n+1$.
\endproclaim
\demo{Proof} Let $r$ be the symmetry rank of $M$. Consider the action of $T^r$
on $M$. By [Sa] the Euler charateristic $\chi(M)>0$. Therefore the
fixed point set of $T^r$ is not empty. Consider the isotropy representation
of $T^r$ at $T_xM\cong \Bbb H^{n}$, where $x$ is a fixed point. Since
$T^r\subset (Sp(n)Sp(1))_x$, so $r\le \text{rank}(Sp(n)Sp(1))=n+1$.
The desired result follows.
\qed\enddemo

\vskip 2mm

\proclaim{Lemma 2.2}
Let $M$ be a positive quaternionic K\"ahler manifold of dimension $4n$ with
symmetry rank $r$. Then there exists a positive quaternionic K\"ahler
submanifold of $M$ of dimension at least $4(r-2)$.
\endproclaim
\demo{Proof} Let $x\in M$ be a fixed point of $T^r$. Consider the isotropy
representation of $T^r$ at $T_xM$. Observe that $\rho _x(T^r)$ has rank
at least $(r-1)$, since the rank of $Sp(1)$ is $1$. Observe that
$\rho _x(T^r)$ acts on $T_xM\cong\Bbb H^n$ quaternionically linearly.
By linear algebra there is a circle subgroup $S^1\subset \rho _x(T^r)$
with fixed point set a quaternionic linear subspace of dimension at least
$4(r-2)$. Let $N$ be the fixed point component of $S^1$ containing $x$.
Clearly the tangent space $T_xN$ is the quaternionic subspace of
$T_xM$ fixed by $S^1$. Since $N$ is totally geodesic
(the fixed point set of an isometric action), for any $y\in N$ there
exists a geodesic $\gamma$ in $N$ joining $x$ and $y=\gamma (1)$. For any
vector
$X\in T_xN$, consider the parallel vector field $X(t)$ along $\gamma$.
By [Ma] Lemma 8.2 we may choose an open set $U$ containing $\gamma$ in $M$
so that the quaternionic structure $I, J,K$ are defined on $U$ and
$IX(t), JX(t), KX(t)$ are all parallel along $\gamma$. Obviously $IX(1),
JX(1), KX(1)$ belongs to $T_yN$. Therefore $T_yN$ is a quaternionic subspace
of $T_yM$. This proves that $N$ is a quaternionic submanifold. By Prop. 1.2
and 1.4 $N$ is a positive quaternionic K\"ahler manifold.
\qed\enddemo

\remark {Remark} It may not be true that any fixed point component of the
above chosen circle subgroup is a quaternionic submanifold. For example,
$\text{Gr}_2(\Bbb C^4)$ admits an isometric circle action with fixed point
set $\text{Gr}_2(\Bbb C^3)\cup \text{Gr}_1(\Bbb C^3)$, the former one is
a quaternionic K\"ahler submanifold, but the latter is not (but K\"ahler).
However, by the proof of Lemma 2.2 we may have a criterion, namely,
if a fixed point component has dimension greater than the middle
dimension $2n$, then it must be a quaternionic submanifold. 
\endremark

An isometric $G$-action on $M$ is called of {\it quaternionic type} if
$\bar \rho _x: G \to Sp(1)_x$ is trivial for any $G$-fixed point $x$.

\vskip2mm

To make our proof more transparent, we now prove Bielawski's result using
our approach, which is essentially included in the following proposition.

\proclaim{Proposition 2.3}
Let $M$ be a positive quaternionic K\"ahler manifold of dimension $4m$ with
symmetry rank $r$. If either $r=m+1$ or $r=m$ and the $T^m$-action is of
quaternionic type. Then $M=\Bbb HP^m$ or $\text{Gr}_2(\Bbb C^{m+2})$.
\endproclaim
\demo{Proof} 
We use induction on the dimension.

By [PS] [HH] we may assume that $m\ge 4$. 
Observe that the fixed point set for the $T^{r}$-action must be isolated.
Let $x\in M$ be a fixed point. By Lemma 2.2 (and its proof) there is a
positive quaternionic K\"ahler $4(m-1)$-submanifold $N^{4(m-1)}$ containing
$x$. Consider the induced $T^r$-action on $N^{4(m-1)}$. Its principal
isotropy
group has rank at most $2$, and has rank $1$ if the $T^r$-action is of
quaternionic type.

Case (i). The principal isotropy group has rank $1$.

There is an effective quotient action of $T^r/S^1=T^{r-1}$ on $N^{4(m-1)}$,
which is quaternionic type if $r=m$. By induction we may assume that
$N^{4(m-1)}=\Bbb HP^{m-1}$ or $\text{Gr}_2(\Bbb C^{m+1})$. By Theorem A3
and Theorem 1.1 the desired result follows.

Case (ii). The principal isotropy group has rank $2$.

The principal isotropy group $T^2\subset T^r$ acts on the normal space of
$T_xN^{4(m-1)}$ in $T_xM$, which may be identified with $\Bbb H$. If
$\bar \rho _x: T^2\to Sp(1)_x$ is trivial, $T^2$ acts on the normal space
quaternionically linearly. For the dimension reason there is a circle
subgroup of $T^2$ acting trivially on the normal space. This is impossible
since the $T^r$-action is effective. In particular, we get that $r=m+1$.

Note that the isotropy representation of the principal isotropy
group does not depend on the choice of the point in $N^{4(m-1)}$. Thus 
$\bar \rho _x: T^2\to Sp(1)_x$ is nontrivial for any $x\in N^{4(m-1)}$.
Note that $\bar \rho _x(T^2)$ acts trivially on $T_x(N^{4(m-1)})$ and
$\bar \rho _x(T^r)=\bar \rho _x(T^2)$. Therefore the
quotient action $T^r/T^2=T^{m-1}$ on $N^{4(m-1)}$ is of quaternionic type.
By induction we may assume that $N^{4(m-1)}=\Bbb HP^{m-1}$ or
$\text{Gr}_2(\Bbb C^{m+1})$. For the same reason as in Case (i) the desired
result follows.
\qed\enddemo

\vskip 4mm

\subhead b). Proof of Theorem B
\endsubhead

\vskip 2mm

We start with several lemmas.

\proclaim{Lemma 2.4}
Let $T^k$ act effectively and quaternionically linearly on a quaternionic
linear space $V\cong \Bbb H^n$ with fixed point set the origin $\{0\}$.
Let $H_1, H_2, \cdots , H_l$ be the set of $(k-1)$-dimensional isotropy
groups, and let $V_1, \cdots , V_l$ be the corresponding subspaces
fixed by the isotropy groups. Let $n_1, \cdots ,n_l$ be the
$\Bbb H$-dimensions of $V_1, \cdots , V_l$. 
Then $l\ge k$, $V=V_1+\cdots +V_l$, and $n=n_1+\cdots +n_l$. 
\endproclaim
\demo{Proof}
By Bredon [Br] III. Theorem 10.12 we know that $V=V_1+\cdots +V_l$ and
$n=n_1+\cdots +n_l$.
It suffices to show that $l\ge k$. Suppose not, the intersection $H_1\cap
\cdots \cap H_l$ is a non-trivial subgroup of $T^k$ which acts trivially on
$V$. A contradiction with the effectiveness of the action.
\qed\enddemo

\vskip2mm

\proclaim{Lemma 2.5}
Let $M$ be a positive quaternionic K\"ahler manifold of dimension $4m$ with
symmetry rank $r$. Assume $b_2(M)=0$. If $r\ge m-2\ge 5$.
Then the $T^r$-action on $M$ has at least an isolated fixed point.
\endproclaim
\demo{Proof}
Since the Euler characteristic $\chi (M)>0$ by [Sa], the fixed point
set of $T^r$ is not empty. Choose a fixed point $x\in M$. Then
$\rho _x: T^r\to Sp(m)$ has an image of rank at least $r-1$.  
Let $N$ be the fixed point component of $\rho _x(T^r)$ containing $x$.
By Lemma 2.2 $N$ is a quaternionic K\"ahler submanifold of dimension $4i$.
Note that $4i+4(r-1)\le 4m$. Thus $i\le 3$.  Clearly we may assume 
$i\ge 1$.

Consider the isotropy representation of $\rho _x(T^r)$ at the normal
vector space of $T_xN\subset T_xM$. By Lemma 2.4 and the assumptions we may
assume
two quaternion subspaces $V_1, V_2$ of dimensions $4$ with codimension one
isotropy groups
$H_1, H_2 \subset \rho _x(T^r)$. Note that $V_1+V_2+T_xN$ is the tangent
space of the fixed point component at $x$ of the action of $H_1\cap H_2$,
saying $N_2$. Note that $N\subset N_2$ has codimension $8$ which extends to
a chain of quaternionic K\"ahler submanifolds
$N\subset N_1\subset N_2$ of codimensions $4$.

\vskip2mm

Case (i). If $i=1$;

Consider the restricted $\rho _x(T^r)$-action on the stratum $N_2$. By the
construction $\rho _x(T^r)$ acts on $T_xN_2$ with principal isotropy group a
codimension $2$ subtorus. Since $\text{dim}N_2=12$, by [HH] we know that
$N_2=\Bbb HP^3$, $\widetilde{\text{Gr}}_4(\Bbb R^7)$ or
$\text{Gr}_2(\Bbb C^5)$.
By Theorem (A3) $b_2(N_2)= b_2(M)=0$. Thus we need only to consider the
former two cases. If $N_2=\Bbb HP^3$ it is obvious. Now we assume that $N_2=
\widetilde{\text{Gr}}_4(\Bbb R^7)$.
By [PS] we know that the quaternionic K\"ahler submanifold $N_1$
must be $\Bbb HP^2$, $G_2/SO(4)$ or $\text{Gr}_2(\Bbb C^4)$.
Since the $\rho _x(T^r)$-action on $\widetilde{\text{Gr}}_4(\Bbb R^7)$ is a
subaction of the standard transitive $SO(7)$-action, one sees easily that any
circle action can not have an $8$-dimensional fixed point component which is
quaternionic. A contradiction.

\vskip2mm

Case (ii). If $i=2$;

Note that $\text{dim }N=8$ and $\text{dim }N_1=12$. Then $\rho _x(T^r)$ acts
on $N_1$ with principal isotropy group rank codimension $1$. As above we
may assume that $N_1=\Bbb HP^3$ or $\widetilde{\text{Gr}}_4(\Bbb R^7)$. For
the same reason as in Case (i) we see that $N_1=\Bbb HP^3$. The desired
result follows.

\vskip2mm

Case (iii). If $i=3$;

By (i) and (ii) we may assume that at every fixed point of $T^r$,
quaternionic K\"ahler submanifolds passing through $x$ fixed by some toric
subgroups, have dimensions at least $12$. Note that $\text{sig}(N)=0$ since
$\text{dim }N=12$ by [HH]. As the Kraines' form $\Omega$ satisfies
$\Omega ^2\ne 0\in H^{8}(N_1)$, we see that $b_8(N_1)\ge 1$. By Theorem 1.3
$\text{sig}(N_1)=
b_8(N_1)>0$.  By [HS] we know that the $T^r$-action has a fixed point
$y\in N_1-N$.

Consider the isotropy representation of $T^r$ on the normal
vector space of $T_yN_1$. Let $H$ be the isotropy subgroup of the linear
action on $T_yN_1$ (i.e. the principal isotropy group of $N_1$). It is easy
to see that $H$ has rank at most $m-3$. Moreover, if $H$ has rank $m-3$,
then $\bar \rho _y: H\subset T^r\to Sp(1)$ is non-trivial and the
restricted $T^r$ action on $N_1$ is quaternionic type.
If $y$ is not an isolated fixed point, we may assume that $T^r/H$ has
rank $1$ and the $T^r/H$-action on $N_1$ is quaternionic type (otherwise
$T^r/H$ may be replaced by a rank one subgroup which is quaternionic type
at $y$). We see that $T^r/H$ has a fixed point component passing through $y$,
which is a quaternionic K\"ahler submanifold in $N_1$ of dimension $\ge 4$.
By Theorem 0.1 this quaternionic K\"ahler
submanifold must intersect with $N$, which is clearly impossible.
By now the proof is complete.
\qed\enddemo

\vskip 2mm

\proclaim{Lemma 2.6}
Let $T^r$ act effectively and isometrically on a quaternionic K\"ahler
manifold $M^{4m}$ with an isolated fixed point $x$. If $r\ge \frac m2+3$,
then there is a chain of quaternionic K\"ahler submanifolds
$$x\subset N^4
\subset N^8\subset N^{12} \subset \cdots N^{4i}
\subset N^{4j}\subset \cdots \subset N^{4s}\subset M^{4m}$$
such that $s\ge r-2$ and the dimensions $4i$, $4j$ of every two nearby
submanifolds in the chain satisfies that $2i-j\ge 2$ for $i\ge 3$.
\endproclaim
\demo{Proof}
Consider the isotropy representation of $\rho _x(T^r)$ at the tangent space
$T_xM\cong \Bbb H^m$. By Lemma 2.4 we get $T_xM=V_1+\cdots +V_l$ where
$l\ge r-1\ge \frac m2+2$ and $V_i$ are all quaternionic subspaces. For the
sake of simplicity, let $V_1, \cdots , V_{i_1}$ be the collection of
subspaces of dimension $4$. By Lemma 2.4 $i_1\ge 4$ since $r-1\ge \frac m2+2$.
Let $V_{i_1+1}, \cdots , V_{i_1+i_2}$ be the collection of subspaces of
dimension $8$, etc.

If $i_1\ge r-2$, consider the sum $V_1+\cdots +V_{r-2}$. Observe  that
$H_1\cap \cdots \cap H_{r-2}\subset \rho _x(T^r)$ has rank at least $1$.
The fixed point component at $x$ of $H_1\cap \cdots \cap H_{r-2}$ is a
quaternionic K\"ahler submanifold of $M$, saying $N^{4r-8}$, of dimension
$4r-8$ whose tangent space at $x$ is $V_1+\cdots +V_{r-2}$.
Clearly the desired chain of submanifolds may be obtained in a similar manner so that
their tangent spaces at $x$ are given by a partial sum $V_1+\cdots +V_{i}$
for $i<r-2$.

If $i_1<r-2$, by Lemma 2.4 $i_2\ge 1$ and $m\ge i_1+2i_2+3i_3+\cdots +
ki_k$, where $4i_k=\text{max} \{\text{dim} V_i: 1\le i\le l\}$. By Lemma
2.4 $i_1+i_2+\cdots +i_k\ge r-1$. If $i_1+i_2
\ge r-2$, choose the quaternionic K\"ahler submanifold
$N^{4i_1+8(r-2-i_1)}$ whose tangent space at $x$ is the sum $V_1+\cdots +
V_{i_1}+V_{i_1+1}\cdots +V_{r-2-i_1}$, which has dimension $>4r-8$. Otherwise
we may proceed this inductively, assuming $i_1+i_2+\cdots +i_{k-1}<r-1$, then
we may take the desired quaternionic K\"ahler submanifold
$N^{4i_1+8(r-2-i_1)}$ whose tangent space at $x$ is the sum $V_1+\cdots +
V_{i_1}+\cdots +V_{i_1+\cdots +i_{k-1}}+\cdots +V_{r-2-i_1-\cdots -i_{k-1}}$.
We may choose the quaternionic K\"ahler submanifolds of $N^{4i_1+8(r-2-i_1)}$
for the chain so that their tangent spaces at $x$ are the partial sums as
above. It is easy to verify that $2i-j\ge 2$ from the construction.
\qed\enddemo

\vskip 2mm

Now we are ready to prove Theorem B.

\vskip 2mm

\demo{Proof of Theorem B}

By Theorem 1.1 we may assume that $b_2(M)=0$.
By Lemma 2.5 we may assume an isolated fixed point $x\in M$ of the
$T^r$-action.

Consider the chain of quaternionic K\"ahler submanifolds defined in
Lemma 2.6. Consider the restricted $T^r$-action on $N^{12}$. 
By the construction the principal isotropy group has rank  $(r-4)$ or
$(r-3)$. In the former case, there is an effective quotient action of
$T^4$ on $N^{12}$. In the latter case we claim that the
$T^r$-action on $N^{12}$ is quaternionic type. Assuming this, by Proposition
2.3 $N^{12}=\Bbb HP^{3}$ or $\text{Gr}_2(\Bbb C^5)$. Therefore
by Theorem (A3) and Theorem 1.1 the desired result follows.

Clearly we need only to consider the case when $\bar \rho _x: T^r\to Sp(1)_x$
is non-trivial.
Since $T^r/T^{r-3}$ acts quaternionically linearly on $T_x (N^{12})\cong
\Bbb H^3$, the restricted representation of the isotropy group $T^{r-3}\subset
T^r@>\bar \rho _x>>Sp(1)_x$ is not trivial. Therefore,  for any
$y\in N^{12}$, the representation of the isotropy group
$T^{r-3}\subset T^r@>\bar \rho _y>>Sp(1)_y$ is not trivial. This implies
readily that the $T^r$-action on $N^{12}$ is quaternionic type.
\qed\enddemo

\vskip 8mm

\head 3. Morse Theory on Path Spaces
\endhead

\vskip4mm

Let $M$ be a complete Riemannian manifold without boundary and let
$f: N\to M\times M$ be an immersed complete submanifold of dimension $n$.

A {\it piecewise smooth path} in $M$ (mod $f$) is a pair $(x, \gamma )$,
where $x\in N$ and a map $\gamma : [0,1]\to M$ such that:

(i) there is a subdivision $0=t_0<t_1<\cdots <t_k=1$
of $[0,1]$ such that each $\gamma |_{[t_{i-1}, t_i]}$ is smooth.

(ii) $(\gamma (0), \gamma (1))=f(x)$.

The set of all piecewise smooth paths (mod $f$) is denoted by $P(M; f)$.

The topology of $P(M;f)$ is taken the induced topology from $N\times P(M)$,
where $P(M)$ is the space of piecewise smooth paths with the metric
topology given by
$$d(\gamma _0, \gamma _1)=\lf\{\int _0^1 (|\dot\gamma_0(t) | -|\dot\gamma_1(t) |
)^2dt\rg\}^{\fr 12} + \max_{0\le t\le 1} d _M(\gamma _0(t), \gamma _1(t)).$$
Note that the integral is well-defined though $\dot\gamma_i$ ($i=1, 2$) may
not be defined at finitely many points in $[0,1]$.

On the space $P(M; f)$ there is an energy function
$E: P(M;f)\to \Bbb R$
given by $$E(x, \gamma )=\int _0^1 | \dot\gamma(t)|^2\,dt.$$
We want to study the topology of $P(M;f)$ using Morse theory for the function
$E$. When $f$ is an embedding, path space with this general boundary condition
was studied in [Gr] using Morse theory on Hilbert manifolds.
Instead, we will use the finite dimensional approximation methods to reduce
to Morse theory in finite dimension.

The tangent space of $P(M;f)$ at $(x, \gamma )$ is defined as the vector space
of piecewise smooth vector fields $W$ along $\gamma$ such that $(W(0), W(1))
\in f_*(T_{x}N)$.

By a standard calculation the first variation of $E$ in the direction $W \in
T_{\gamma }P(M,f)$, is given by
$$ \al
\frac 12 E_*(W)=&-\int_0^1\lf<W(t),\fr{D}{dt}\dot\g(t)\rg>\,dt+\cr
&\lf<W(1),\dot
\g(1)\rg>-\lf<W(0),\dot\g(0)\rg>+\sum_{i=1}^{k-1}\lf<W(t_i),\dot\g_-(t_i)-\dot\g_+(t_i)\rg>,\eda
$$
where $\dot\g_-$ is the left derivative and $\dot\g_+$ is the right
derivative of
$\g$, and $\fr{D}{dt}$ is the covariant derivative along $\g$.
So if $(x, \gamma )$ is a critical point of
 $E$ then:

\noindent (i) $\gamma $ is a smooth geodesic;

\noindent (ii) $(\dot \gamma(0),-\dot \gamma(1))$ is perpendicular to
$f_*(T_x(N))$ at $f(x)=(\gamma(0),\gamma(1))$.

Let $W_1,W_2\in T_\gamma P(M;f)$. If $\gamma$ is a critical point of $E$ we
consider any variation $h(t,s,u)$ of $\g$ with $\fr{\p h}{\p s}(t,0,0)=W_1(t),
\fr{\p h}{\p u}(t,0,0)=W_2(t)$. Then the
second variation of $E$ along $\gamma$, denoted by $E_{**}(W_1, W_2)$, is as
follows:
  $$\al
\frac 12 E_{**}(W_1, W_2)=&\int_0^1\lf\{\lf<\fr{DW_1}{dt},\fr{DW_2}{dt}\rg>
-\lf<R(\dot\g,W_1)W_2,\dot\g\rg>\rg\}+\cr
&\lf<\fr{D}{ds}\fr{\p h}{du}(1,0,0),\dot\g(1)\rg>-
\lf<\fr{D}{ds}\fr{\p h}{du}(0,0,0),\dot\g(0)\rg>\cr
=&\int_0^1\lf\{\lf<\fr{DW_1}{dt},\fr{DW_2}{dt}\rg>
-\lf<R(\dot\g,W_1)W_2,\dot\g\rg>\rg\}+\cr
&\Bigl<\a\bl((W_1(0),W_1(1)),(W_2(0),W_2(1))\br),\bl(-\dot\g(0),\dot\g(1)\br)\Bigr>,\eda
$$
where $\a$ is the second fundamental form of the immersion  $f:N\to M\times M$.

Let $P_c(M;f)=E^{-1}([0,c)) \subset P(M;f)$.
Following Milnor-Morse
we define a finite dimensional approximation to $P_c(M;f)$ as follows:

Choose some subdivision $0=t_0<t_1<\cdots < t_k=1$ of $[0,1]$. Let $B$
 be the
subspace of ${P_c}(M;f)$ such that

(i) $f(x)= (\gamma (0), \gamma (1))$

(ii) $\gamma |_{[t_{i-1}, t_i]}$ is a geodesic for each $i=1, \cdots , k$.

\vskip2mm

\proclaim{Theorem 3.1}
Let $M$ be a complete Riemannian manifold, and $f: N\to M\times M$ be
an immersion where $N$ is a closed Riemannian manifold. Let $c$ be a fixed
positive number  such that $P_c(M; f)$ is not empty. Then for all
sufficiently fine subdivision $0=t_0<t_1<\cdots < t_k=1$ of $[0,1]$ the set
$B$ can be given the structure of a smooth
finite dimensional manifold.
\endproclaim
\demo{Proof}
[Mi], Sect. 16.
\qed\enddemo

\vskip2mm

Let $E|_B: B\to \Bbb R$ be
the restriction of the energy function $E$.

\vskip 2mm

\proclaim{Theorem 3.2} Let $f,N,M$ be as in Theorem 2.1. Then
$E|_B$ is a smooth function. For each $a<c$ the set
$(E|_B)^{-1}([0,a])$ is compact, and
 $(E|_B)^{-1}([0,a))$ is a deformation retract
of the set $P_a (M; f)$. The critical points of $E|_B$ are precisely the same as
the
critical points of $E$ in $P_c(M; f)$, that is, the pairs $(x, \gamma
)$, where $\gamma$  is a smooth geodesic in $M$ such that $(\dot \gamma (0),
-\dot \gamma (1))$ is normal to $f_*(T_x(N))$ and the energy is less than $c$.
The Hessians of $E|_B$ and $E|_{{P_c}(M;f)}$ have the same
index at each critical point $(x, \gamma )$.
\endproclaim
\demo{Proof}
[Mi], Sect. 14 and Sect. 16.
\qed\enddemo

\vskip 2mm

\proclaim{Lemma 3.3}
Suppose that every nontrivial critical point $(x, \gamma )$ of $E$ has
positive index. Then $f^{-1}(\Delta )$ is not empty.
\endproclaim
\demo{Proof}
Suppose $f^{-1}(\Delta )=\emptyset$. Writing $f(x)=(f_1(x),f_2(x))$,
from the compactness
of $N$ we conclude that there exists $x\in N$ such that $
d(f_1(x),f_2(x))=\de$ is a positive minimum.
In particular $P(M; f)$ contains no constant path  and the energy
function on $P(M;f)$ assumes a  positive minimum $\de^2$ at $(x,\g)$,
where $\g$ is a minimal geodesic
joining $f_1(x)$ and $f_2(x)$.
For this minimal critical point $(x,\g)$ the index is clearly zero.
A contradiction to the assumption.
\qed\enddemo

\vskip 2mm

We also need the following lemma of Milnor for finite dimensional manifolds.

\vskip 2mm

\proclaim{Lemma 3.4}
Let $X$ be a finite dimensional smooth manifold and $f: X\to \Bbb R$ be a real
function
with minimal value zero. Suppose that for any $a$ the sublevel set $X_{\le
a}=f^{-1}([0, a])$
is compact. Assume that the set $X_0$ of minimal points has a neighborhood $U$
with
a deformation retraction $r: U\to X_0$, and that all nontrivial critical points
have
indices greater than $\lambda _0 \ge 0$. Then $X$ has the homotopy-type
of a CW-complex
by attaching cells of dimensions at least $\lambda _0+1$ to $X_0$. In
particular,
$\pi _j(X, X_0)=0$ for $0\le j\le \lambda _0$.
\endproclaim
\demo{Proof}
[Mi], Sect. 22.
\qed\enddemo

\vskip 2mm

\proclaim{Theorem 3.5}
Let $M, N_i$, $i=1, 2$, be compact quaternionic K\"ahler manifolds and
$\Delta$ be the diagonal of $M\times M$. Let $f_i: N_i\to M$ be quaternionic
immersions. Let $f=(f_1, f_2)$.  Let $P_0=f ^{-1}(\Delta )$.  If\ \,every nontrivial critical point of $E$ on
$P(M;f)$ has index $\lambda >\lambda _0\ge 0$, then $P(M;f)$ has the
homotopy-type of a CW-complex obtained by attaching  cells of dimensions at
least $\lambda _0+1$ to $P_0$. In particular,  the relative homotopy groups
$\pi _j(P(M; f), P_0)=0$ for $0\le j\le \lambda _0$.  \endproclaim
\demo{Proof}
Note that each point $x$ in $P_0$ could be associated with
a constant path at $p_1f(x)$ ($=p_2f(x)$). So $P_0$
can be identified with the constant geodesics in $P(M; f)$.
It suffices to prove that for any large value $c$, $P_c(M, f)$ has
the homotopy-type of a CW-complex obtained by attaching
cells of dimensions at least $\lambda _0+1$ to $P_0$. By Theorem 3.2
$P_c (M, f)$ deformation retracts to $B$.
Moreover, the index for every nontrivial critical point for the restricted
energy function
is greater than $\lambda _0$ too. The space $P_0$ is clearly inside
$B$. To apply Lemma 3.4 it suffices to prove that there is a neighborhood $U
\subset B$ of $P_0$ and a retraction $r: U \to P_0$. 

Since a quaternionic immersion has to totally geodesic (cf. Prop. 1.4),
the minimal set $P_0$ (resp. $B$) may be identified with the submanifold
$\{(x, f(x), \cdots ,f(x)): x\in f^ {-1}(\Delta) \}$ (resp. an open
submanifold of $N\times (M\times M) \times \cdots \times (M\times M)$).
Therefore an open regular neighborhood of this submanifold is a
desired open neighborhood $U$. The desired result follows.
\qed\enddemo

\vskip2mm

\head 4. Indices Theorems For critical points of Energy Functions
\endhead

\vskip2mm

Let $V$ be a $4l$-dimensional quaternionic linear space and let $Q$ be a real
symmetric bilinear form on $V$. Recall that the {\it index } of $Q$ is
defined by the dimension of a maximal linear subspace  $U\subset V$ so that
$Q|_U$ is negative definite.

\proclaim{Lemma 4.1} Let $V$ be a $4l$-dimensional quaternionic linear space
with invariant inner product $\lf <.,.\rg>$ (i.e. $I, J, K=IJ$ preserve the
inner product.) Let $Q(.,.)$ be a real symmetric bilinear form on $V$.
Assume that for all nonzero vector $X\in V$ it holds that
$$
Q(X,X)+Q(IX, IX)+Q(JX, JX)+Q(KX, KX)<0
$$
Then the index $\la$ of $Q$ satisfies $\la\ge l$.
\endproclaim
\dem Let $A: V\to V$ be the symmetric linear map uniquely determined by
$\lf<Av,w\rg>=Q(v,w)$ for $v, w\in V$. Let $\Cal B$ be an orthonormal basis
of eigenvectors of the linear map $A:V\to V$. Consider the linear subspace
$W$ spanned by the vectors $v\in \Cal B$ with $Q(v,v)\ge 0$. Observe that
the index of $Q$ is equal to $4l-\text{dim} (W)$. Therefore it suffices to
prove that $\text{dim} (W)\le 3l$. Suppose not, then the intersection
$W\cap IW \cap JW\cap KW$ is a nontrivial linear subspace. For any
nontrivial $v_0=Iv_1=Jv_2=Kv_3\in W\cap IW \cap JW\cap KW$, where
$v_i\in W$, we have $Q(v_i, v_i)\ge 0$. Therefore
$Q(v_0, v_0)+Q(Iv_0, Iv_0)+Q(Jv_0, Jv_0)+Q(Kv_0,Kv_0)\ge 0$.
 A contradiction. The desired result follows.
\qed\edo

Let $M$, $N_i$, $i=1, 2$, be compact quaternionic K\"ahler $4m$-manifolds.
Let $N=N_1\times N_2$, and let $f_i: N_i\to M$ be quaternionic immersions.
Set $f=(f_1, f_2)$ and $4n=4n_1+4n_2$ for the real dimension of $N$. 
Let $W_1, W_2\in T_{(x, \gamma )}P(M;f)$ be tangent vectors at $(x, \gamma )$.
Let $\a$ be the second fundamental form of $f$ in  $M\times M$. Recall that
the second variation of the energy function $E$ along a critical
point $(x,\gamma)$  reads
$$\al
\frac 12 E_{**}(W_1, W_2)
=&\int_0^1\lf\{\lf<\fr{DW_1}{dt},\fr{DW_2}{dt}\rg>
-\lf<R(\dot\g,W_1)W_2,\dot\g\rg>\rg\}+\cr
&\Bigl<\a\bl((W_1(0),W_1(1)),(W_2(0),W_2(1))\br),\bl(-\dot\g(0),\dot\g(1)\br)\Bigr>.
\eda\tag 4.2
$$
Now let $W$ be parallel along $\gamma$ with $\Cal W=(W(0),W(1))\in f_*(T_xN)$.
By [Ma] Lemma 8.2 we may assume that $I\Cal W$, $J\Cal W$ and $K\Cal W$ are
all parallel and tangent to $N$ at $x$, since $f_i$ are quaternionic
immersions. Set
$$
\eta=\bl(-\dot\g(0), \dot\g(1)\br).$$
Since $(x,\g)$ is a critical point, $(\dot \g (0), -\dot \g (1))$ is
orthogonal to $\Cal W$, $I\Cal W$, $J\Cal W$ and $K\Cal W$. By the formula
(4.2) we get
\vskip 4mm

$
\frac 12 E_{**}(W, W)+ \frac 12 E_{**}(IW, IW) +\frac 12 E_{**}(JW, JW)
+\frac 12 E_{**}(KW, KW)  = $

$=\!\int_0^1-\left\{\lf<R(\dot\g,W)W,\dot\g\rg>\! +\!\lf<R(\dot\g,IW)IW,
\dot\g\rg>\! + \!\lf<R(\dot\g,JW)JW,\dot\g\rg>\!+ \!\lf<R(\dot\g,KW)KW,
\dot\g\rg>\! \right\}
$

\vskip 4mm

Let $\Cal V$ be the real linear space spanned by vectors $(W(0),W(1))$,
where $W$ is a parallel vector field  along $\g$ so that $W$ orthogonal to
$\dot\g$, $I\dot\g$, $J\dot\g$ and $K\dot\g$.
The quaternionic dimension of $\Cal V$ is $m-1$. If $\Cal W=(W(0),W(1))\in
\Cal V$ then $I\Cal W$, $J\Cal W$ and $K\Cal W\in \Cal V$ and
$\Cal W, I\Cal W$,
$J\Cal W$, and $K\Cal W$ are all orthogonal to $\eta=(-\dot\g(0),\dot\g(1))$,
$I\eta$, $J\eta$ and $K\eta$. Therefore the quaternionic dimension
$$\dim _\Bbb H\bl(\Cal V\cap f_*(T_xN)\br)\ge m-1+n-(2m-1)=n-m.\tag 4.3
$$
Recall that if $\Cal W\in \Cal V\cap f_*(T_xN)$ then $W$ is a tangent vector
of $P(M, f)$ at ${(x,\g)}$.

\proclaim{Theorem 4.2}
Let $M$ be a positive quaternionic K\"ahler manifold of dimension $4m$. Let
$f_i: N_i\to M$, $i=1, 2$, be quaternionic immersions of closed quaternionic
K\"ahler manifolds $N_i$ of dimensions $4n_1$ and $4n_2$. Set $n=n_1+n_2$.
Let $f=(f_1, f_2)$ and $N=N_1\times N_2$. If $(x,\gamma)$ is a nontrivial
critical point for the energy function $E$ on $P(M, f)$, then the index
$\la$ of $E_{**}$ at $(x,\g)$ satisfies $$\la \ge n-m+1.$$
\endproclaim
\dem
Let $V\subset T_{(x,\g)}P(M,f)$ be the quaternion
vector space spanned by parallel vector fields $W$ along $\g$ orthogonal
to the quaternion line $\dot\g\wedge I\dot \g \wedge J\dot\g\wedge K\dot\g$
such that $(W(0),W(1))\in \Cal V\cap f_*(T_xN)$. Note
that $V$ and $\Cal V\cap f_*(T_xN)$ are $\Bbb H$-isomorphic, hence by (4.3) we
get $\dim_{\Bbb H}(V)\ge n-m$.

We claim that in our case $\dim_{\Bbb H}(V)\ge n-m+1$.

Indeed, if
$x=(x_1,x_2)$, it is easy to see that the critical point $(x,\g)$ satisfies
$\dot \gamma (0)\perp (f_1)_*(T_{x_1}N_1)$ and $\dot \gamma (1)
\perp (f_2)_*(T_{x_2}N_2)$. Therefore $(\dot \gamma (0),0)$, $(I\dot \gamma
(0), 0)$, $(J\dot \gamma (0),0)$ and $(K\dot \gamma (0),0)$ are all
orthogonal to $(f_1)_*(T_{x_1}N_1)\times (f_2)_*(T_{x_2}N_2)=f_*(T_xN)$;
similarly $(0, \dot \gamma (1))$, $(0, I\dot \gamma (1))$,
$(0, J\dot \gamma (1))$ and $(0, K\dot \gamma (1))$ are all orthogonal to
$f_*(T_{x}N)$. Obviously, all the vectors $(\dot \gamma (0), 0)$,
$(I\dot \gamma (0), 0)$, $(J\dot \gamma (0), 0)$, $(K\dot \gamma (0), 0)$, and
$(0, \dot \gamma (1))$, $(0, I\dot \gamma (1))$, $(0, J\dot \gamma (1))$
and $(0, K\dot \gamma (1))$ are orthogonal to $\Cal V$. Therefore
$$\dim_{\Bbb H}(V)=\dim_{\Bbb H}(\Cal V\cap f_*(T_{x}N)) \ge n-m+1$$

Endow $V$ with the inner product $\lf<X,Y\rg>:=\lf<X(0),Y(0)\rg>$. Consider
the real symmetric bilinear form
$Q(v,w)=E_{**}(v,w)$.
Note that the index of $E_{**}$ is not less than the index of $E_{**}$
restricted to $V$. Therefore by Lemma 4.1 it suffices to prove that
for any nontrivial vector $X\in V$ it holds that
$$Q(X,X)+Q(IX, IX)+Q(JX, JX)+Q(KX, KX)<0,\tag 4.4 $$
where $I, J, K$ defines the quaternionic structure on $V$.

By Proposition 1.2 the left side of (4.4) equals to
$$
-\!\int_0^1\!\! \frac \mu {m+2}\|X\|^2 \|\dot \g\|^2  <0
$$
where $\mu$ is the Einstein constant of $M$.
The desired result follows.
\qed\edo

\vskip 5mm

\head 5. Proofs of Theorem A
\endhead

Let $M$ be a positive quaternionic K\"ahler manifold. Let $P(M)$ denote the
space
of $C^ k$ paths $\gamma : [0,1]\to M$. Let $\pi : P(M)\to M\times M$ be
the projection to the pair of end points (i.e. $(\gamma (0), \gamma (1))$).
It is a standard result in topology that $\pi$ is a Serre fibration
$$\Omega M\to P(M)@> \pi >> M\times M$$
with fiber the loop space $\Omega M$ with a fixed base point.
By definition we know that the projection to the first factor
$p: P(M, f)  \to N$ is exactly the pullback fibration by
$f: N\to M\times M$ from $\pi : P(M)\to M\times M$.
In particular, there is a homotopy exact sequence
$$\cdots \to \pi _i(P(M; f))\to \pi _i(N)\to \pi _{i-1}(\Omega M)
\to \pi_{i-1}(P(M;f))\to \cdots \tag 5.1$$

\vskip6mm

\demo{Proof of Theorem A}

By Lemma 3.3 and Theorem 4.2, (A1) follows.

Note that $\pi _{i-1}(\Omega M)=\pi _i(M)$. In particular, $\Omega M$
is path connected since $M$ is simply connected (Theorem 1.1). Then
$P(M;f)$ is path
connected. If $n>m+1$, by Theorem 3.5 and Theorem 4.2,  $P(M; f)$
has the homotopy-type of a CW-complex, obtained from $f^{-1}(\Delta )$ by
attaching cells of dimensions at least $2$. Thus $f^{-1}(\Delta )$ is path
connected. This proves (A2).

Let $p _1: P(M; f)\to N_1$ denote the composition of the bundle projection
$p$ with the projection from $N_1\times N_2\to N_1$  to the first factor.
Note that
$p_1$ is a fibration with fiber $V$, the pullback fibration
fitting in the commutative diagram below
$$\CD
\Omega M& & =& &\Omega M \\
@VVV & &  @VVV\\
V & @>>> & P(M,*)  \\
@V{p_2} VV & &  @V \pi VV\\
N_2& @>f_2>> & M
\endCD
$$
where the right side is the principal path fibration. Note that
$P(M,*)$ is contractible. By the above diagram $V$ is the homotopy fiber
of the map $f_2: N_2\to M$. Therefore by the fibration homotopy exact
sequence it follows that
$$\pi _i(V)\approx \pi _{i+1}(M, N_2)$$
for all $i$. This together with the homotopy exact sequence for the
fibration
$$\CD
V \to P(M; f) \\
@V{p_1}VV \\
N_1
\endCD
$$
and the long exact sequence for the map $i_1: f^ {-1}(\Delta )=N_1\cap N_2
\to N_1$ gives a commutative diagram:
$$
\CD
\pi _{i+1}(N_1) & \to & & \pi _{i+1} (N_1, f^{-1}(\Delta )) & & \to &
\pi _{i}(f^{-1}(\Delta )) & & \to & \pi _{i}(N_1) & & \to &
\pi _{i}(N_1,  f^{-1}(\Delta )) \\
@V{=}VV & @V{(f_1)_*}VV &  @V{surj}VV &  @V{=}VV & @V{(f_1)_*}VV\\
\pi _{i+1}(N_1) & \to & & \pi _{i+1} (M, N_2) & & \to & \pi _{i}(P(M; f))
& & \to & \pi _{i}(N_1) & & \to & \pi _{i}(M, N_2)
\endCD
$$
The middle homomorphism is surjective for $i\le n-m$ since
by Theorem 3.5 and Theorem 4.2 $\pi _i(P(M;f), f  ^{-1}(\Delta ))=0 $.
From the 5-lemma the commutative diagram above implies that
$$(f_1)_*: \pi _i(N_1, f^{-1}(\Delta ))\to \pi _i(M, N_2)$$
is an isomorphism for all $i\le n-m$ and a surjection for $i=n-m+1$.
The desired result follows.
\qed\enddemo

\vskip4mm

\demo{Proof of Theorem 0.2}

Since $f: N\to M$ is an immersion, it suffices to show that $f$ is
a one-to-one map. Note that $(f,f)^{-1}(\Delta)=\{(x, x), x\in N\}
\cup \{(x, y): f(x)=f(y), x\ne y)\}$. Hence, if $f$ is not one-to-one,
then $f^{-1}(\Delta)$ is not connected; a contradiction to (A2).
\qed\enddemo

\enddocument